\RequirePackage{fix-cm}

\documentclass[smallextended]{svjour3}       % onecolumn (second format)
\smartqed  % flush right qed marks, e.g. at end of proof

\usepackage{graphicx}
\usepackage{epstopdf}

\usepackage{amsmath}
\usepackage{mathrsfs}

\usepackage{algorithm}
\usepackage{algorithmic}

\usepackage{booktabs}
\usepackage{enumerate}

\begin{document}
%
% paper title
% Titles are generally capitalized except for words such as a, an, and, as,
% at, but, by, for, in, nor, of, on, or, the, to and up, which are usually
% not capitalized unless they are the first or last word of the title.
% Linebreaks \\ can be used within to get better formatting as desired.
% Do not put math or special symbols in the title.
\title{A robust optimization and an area decomposition approach to
	large-scale network planning under uncertainty
}

%% To specify the authors when (number of affiliations <= 2)

\author{Maynara A. Aredes,     \and
        Tiago Andrade,         \and
        Gerson C. Oliveira,    \and
        Maria L. Latorre,      \and
        Daniela Bayma,         \and
        Silvio Binato
}

%\authorrunning{Short form of author list} % if too long for running head

\institute{All authors are with PSR at:\at
              Praia de Botafogo, 370 - Botafogo, Rio de Janeiro - RJ, 22250-040, Brazil \\
              Tel.: +55-21-3906-2100\\
              %Fax: +123-45-678910\\
              \email{\{maynara,tiago.andrade,gerson,lujan,daniela,silvio\}@psr-inc.com}           %  \\
                                          \\
              %\\
              %\hfill
              %This work was presented in ICSP 2016 - Buzios, Brazil
%             \emph{Present address:} of F. Author  %  if needed
%           \and
%           S. Author \at
%              second address
}

\date{April 2, 2020}

% make the title area
\maketitle

% As a general rule, do not put math, special symbols or citations
% in the abstract
\begin{abstract}
This paper describes an area decomposition and a robust optimization strategy for the planning of large-scale transmission networks considering uncertainty in variable renewable energy (VRE) and hydropower generation. The strategy is illustrated with a study of Brazil’s power grid expansion from 2027 to 2036. The country’s network has over 10 thousand buses and 14 thousand circuits and covers an area equivalent to the USA. A clustering scheme based on bus marginal costs is used to decompose the grid into seven regions. Each regional planning is then solved by a heuristic and Benders decomposition methods to produce a robust least-cost plan for 3,000 dispatch scenarios.
\end{abstract}

\keywords{
Heuristic \and Robust optimization \and Transmission expansion problem
}

\section{Introduction}
\label{sec: intro}

Planning the Brazilian transmission system is a challenging analytical and computational task for several reasons: (i) size: the current HV grid has 10 thousand buses, 14 thousand circuits and spans an area like that of the USA or EU; (ii) investment cost: about 100 thousand km of HV lines were built in the past twenty years, and 50 thousand km more are planned for the next seven years; (iii) operating scenarios: hydropower corresponds to 70\% of Brazil’s 160 GW installed capacity. In order to exploit the country’s rainfall diversity, hydro plants were built in several different river basins. Consequently, there are many possible power transfer scenarios between “wet” and “dry” regions, which must be considered in the planning process.

With the recent insertion of VRE resources (20 GW of wind and 4 GW of solar generation) into the system (with much more being envisioned), the power flow scenarios to be considered in grid planning have become even more diverse. The main reasons are: (a) the synergy between wind power and inflows in the country’s Northeast region (wind production is highest during the dry season) affected the import/export patterns, and new congestion zones have appeared; (b) hydro reservoirs are used to manage renewable variability (“water batteries”).

This paper describes the new methodologies and analytical tools used to determine the least-cost expansion of Brazil’s transmission network, considering the above challenges. The first step in the planning process is to use a clustering scheme on the bus marginal costs to divide the country into a given number of areas, \cite{rsampaio2019community}. Then, the reduced system configuration is used in the regional planning process to determine the interconnections reinforcements (i.e. only tie-line overloads among those areas are considered; all other circuit limits are relaxed). 

This regional planning  problem is formulated as the co-optimization of generation and area interconnection investment costs, plus expected operation costs over many hydrological and renewable scenarios. The problem is solved by stochastic Benders decomposition, where the master is a mixed integer (MIP) investment module that determines candidate generation and interchange capacities. In turn, the subproblem is a stochastic multistage operation module that sends optimality cuts to the investment module (see \cite{pereira1991multi,campodonico2003expansion,fern2019stochastic} for a detailed description of this regional planning scheme).

The main focus of this paper is on the next step: given the generation and area interconnection reinforcements, determine the optimal grid expansion plan for each area. This problem is formulated as the minimization of investment costs (circuit and transformers) required to avoid overloads for a set of operating scenarios (robust optimization). Each scenario comprises two vectors, one with bus generation injections (renewables, hydro and thermal plants) and bus load injections; and the other with the area boundary injections (power flows in the tie-lines). The scenarios are derived from the stochastic operational simulation of the Brazilian system (including generation and area interconnection reinforcements) \cite{campodonico2003expansion}. The study’s simulation covered the years 2027 to 2036 and had 84 inflow/VRE scenarios, 12 months, and 3 load blocks per month, resulting in 3024 operation scenarios for each area per year.

Because of Brazil’s network size and high load growth, it is not necessary to have a multi-year assessment of transmission reinforcements. Therefore, the planning procedure is applied sequentially for each year of the study period. The following heuristic procedure initially obtains the expansion for a given year: (i) rank the 3024 scenarios based on both severity and geographical diversity of overloads; (ii) solve the large-scale MILP (minimize investment costs subject to power flow equations and constraints) for the K “worst” scenarios (typically, K=5); (iii) implement the resulting reinforcements; (iv) re-rank the scenarios and return to steps (ii)-(iv) until all overloads are eliminated; (v) because of the “myopic” nature of the heuristic procedure, some circuit investments may be redundant in the final solution; another heuristic eliminates those redundancies: rank the reinforcements by decreasing costs and eliminate them if their simulated removal does not lead to overloads in any scenario.
The above heuristic procedure often produces a nearly optimal transmission plan, which can be used in case the Benders decomposition algorithm, applied next, does not converge in the allotted time. The Benders algorithm has a MIP investment module that produces candidate plans; and 3024 operational modules, each with a linearized OPF that sends a feasibility cut if there are overloads for the operational scenario (see \cite{latorre2019stochasticrobust} for a detailed description). The convergence of the Benders is accelerated by incorporating the feasibility cuts produced by the MILP problems of the heuristic procedure.

Section~\ref{sec: prob. and met.} describes the transmission network expansion problem and reviews solution methodologies. Section~\ref{sec: heuristic} describes the heuristic scheme for the robust transmission planning of each area. Section~\ref{sec: case study} presents the Brazilian study assumptions, computational performance, and results. Finally, Section~\ref{sec: conclusions} discusses the conclusions.

\section{Solution Methodology}
\label{sec: prob. and met.}

Transmission expansion planning (TEP) is a topic of great interest, first studied by Garver~\cite{garver1970transmission} in 1970. Since then, several methodologies have been proposed. Extending the static model of Garver, a dynamic model was introduced by Dudu and Merlin~\cite{dodu1981dynamic}. TEP models are classified as static when they aim to solve a single-year expansion~\cite{chen2014robust,ruiz2015robust,dehghan2015robust,minguez2016robust}, and as dynamic when solve a multi-year problem \cite{chen2016robust,garcia2016dynamic}. As seen, we solve a multi-year problem by solving static chronological problems, fixing the previous investments.

Another classification for TEP models is whether they consider uncertainty or not. Uncertainty is usually dealt based on stochastic programming (SP) \cite{park2015stochastic,park2013transmission,liu2017multistage,conejo2017long}, and robust optimization (RO) \cite{jabr2013robust,chen2014robust,ruiz2015robust,minguez2016robust,chen2016robust,dehghan2017adaptive}. As seen, we used a SP approach in the generation and area interconnection planning, and a RO approach for the transmission planning of each area. 

The TEP is nonconvex and hard to solve, even if the linear approximation is used as is the usual case since the problem is a large-scale combinatorial problem. Therefore, along with exact methods \cite{latorre2019stochasticrobust,orfanos2012transmission,de2008transmission,zhao2009flexible,Zhang2017Robust,Zhan2017Fast,binato2001new}, many heuristics \cite{chanda1998reliability,sousa2011combined,pereira1985application,romero2003analysis}, and metaheuristics  \cite{da2011performance,leou2011multi,liu2001application,binato2001greedy,da2010reliability} were proposed to solve it. In this paper, we present a greedy heuristic to obtain high-quality solutions in feasible computational time.

%In this paper, we present a greedy heuristic that first do constructive phase to obtain a feasible investment plan, and then a local search phase to eliminate redundant investments. We use the heuristic to solve a large-scale robust model for TEP and compare it with an exact method based on Benders decomposition.

Some papers related to the current work include \cite{bahiense2001mixed}, where a MIP disjunctive model was proposed; \cite{binato2001greedy}, which presented a greedy randomized adaptive search procedure (GRASP); \cite{binato2001new,granville1988mathematical}, which proposed a Benders decomposition for TEP model and showed how to find adequate values for the big M in the disjunctive constraints; and \cite{latorre2019stochasticrobust}, which presented the Benders decomposition for the transmission expansion, without the heuristic to accelerate it. For a more detailed TEP literature review, the reader should refer to Hermmanti et al. \cite{hemmati2013comprehensive}, and Mahdavi et al.~\cite{mahdavi2018transmission}.

\subsection{Mathematical formulation}
\label{sec: formulation}

The robust TEP for a given year can be formulated as:

\begin{align}
\min & \ \sum_j c_j x_j \\
\text{s.t.} & \ z^s(x) \leq 0 \quad \forall s = 1,\dots,S
\end{align}

Where $j=1,\dots,J$ indexes the reinforcement candidates; $c_j$ is the investment cost of candidate $j$; and $x_j$ is the binary investment decision.

The functions $z^s (x)$, $s = 1,\dots,S$ indicate the feasibility of a given expansion plan $x= \{x_j\}$ for each operating scenario $s = 1,\dots,S$. This feasibility is represented as the solution of the following linear programming (LP) problem:

\begin{align}\label{eq1}
z^s (x) = \min & \ e' r^s \\
\text{s.t.} & \ [I] (f_e^s + f_c^s) + r^s = d^s - g^s \\
& \ f_e^s = \Gamma_e [I]' \Theta^s \\
& \ |f_c^s - \Gamma_c [I]' \Theta^s| \leq M_c (1-x) \\
& \ - \bar f_e \leq f_e^s \leq \bar f_e \\
& \ - \bar f_c x \leq f_c^s \leq \bar f_c x \\
& \ r^s \leq d^s\label{eq2}
\end{align}

Table~\ref{tab: Notation} summarizes the notation for (3)-(9) . 

\begin{table}[!ht]
	% increase table row spacing, adjust to taste
	\renewcommand{\arraystretch}{1.3}
	\centering
	\caption{Notation}
	\label{tab: Notation}
	\begin{tabular}{c c}
		\toprule
			$[I]$ & Branch-node incidence matrix $[I]'$ is the transpose) \\
			$r^s$ & Bus load curtailment for scenario $s$ \\
			$d^s$ & Bus load, scenario $s$ \\
			$g^s$ & Bus generation, scenario $s$ \\
			$f_e^s$ & Flow for existing circuits, scenario $s$ \\
			$f_c^s$ & Flow for candidate circuits, scenario $s$ \\
			$\Gamma_e$ & Suceptance matrix for existing circuits \\
			$\Gamma_c$ & Suceptance matrix for candidate circuits \\
			$\Theta^s$ & Bus voltage angle, scenario $s$ \\
			$\bar f_e$ & Flow limit for existing circuits \\
			$\bar f_c$ & Flow limit for candidate circuits \\
			$M_c$ & “big-M” coefficient for candidates \\	
		\bottomrule
	\end{tabular}
\end{table}

The objective of $z^s (x)$ is to minimize the load curtailment $r^s$ required to eliminate circuit overloads. As a consequence, $z^s (x)=0 \Leftrightarrow $ feasibility (no overloads). The first equation represents the power balance in each bus (Kirchhoff’s 1st law); the second equation is Kirchhoff’s 2nd law for the existing circuits; the third set of constraints represents the 2nd law for candidates (disjunctive formulation, with a customized “big M” \cite{binato2001new}). The last constraints are bounds for circuit flow and bus load curtailment.

The TEP problem can be solved directly as a large-scale MIP by incorporating the constraints of the minimum load curtailment LP problem into the TEP formulation; or iteratively, by a Benders decomposition scheme, where the constraints $z^s (x) \leq 0$ are approximated as feasibility cuts.

The above limitations motivated the scheme presented in this paper to solve MIP problems sequentially, with a relatively small number of “severe” scenarios in each step, thus producing a feasible solution in a relatively short time; and then apply the Benders decomposition scheme for the remainder of the allowed time, using feasibility cuts produced during the heuristic process to accelerate convergence. As mentioned, experience has shown that the feasible solutions produced by the heuristic alone are often nearly optimal, and this will be shown in the study case.

\subsection{Heuristic Procedure}
\label{sec: heuristic}

This procedure is implemented in the following steps.

\begin{enumerate}[(i)]
	\setlength\itemsep{1em}
	\item
	rank the $S$ scenarios by decreasing severity (the ranking criterion considers the amount and geographical diversity of overloads);  % put exact criterion
	
	\item \label{item: heuristic severe investment}
	select the $k = 1,\dots, K$ most severe scenarios in the ranking; solve the robust TEP problem for this subset (as mentioned, typically, $K=5$);
	
\begin{align}
	\min & \ \sum_j c_j x_j \\
	\text{s.t.} & \ z^k \leq 0 \quad \forall k = 1,\dots,K
\end{align}
	
	\item
	incorporate to the network the reinforcements resulting from the optimal solution of step~(\ref{item: heuristic severe investment});	
	
	\item
	re-rank the scenarios S and return to step~(\ref{item: heuristic severe investment}) until all overloads are eliminated;
	
	\item
	
	eliminate redundancies: rank the reinforcements by decreasing costs and eliminate them if their simulated removal does not lead to overloads in any scenario.
\end{enumerate}

\section{Case study: Brazilian power system}
\label{sec: case study}

The methodologies described in the previous section were applied to the expansion planning of Brazil’s network from 2027 to 2036. The reason for starting in 2027 is that the generation and transmission reinforcements until 2026 are already under construction or have been contracted by past generation and transmission auctions. 

\begin{figure}[!ht]
\centering
\includegraphics[width=2.7in]{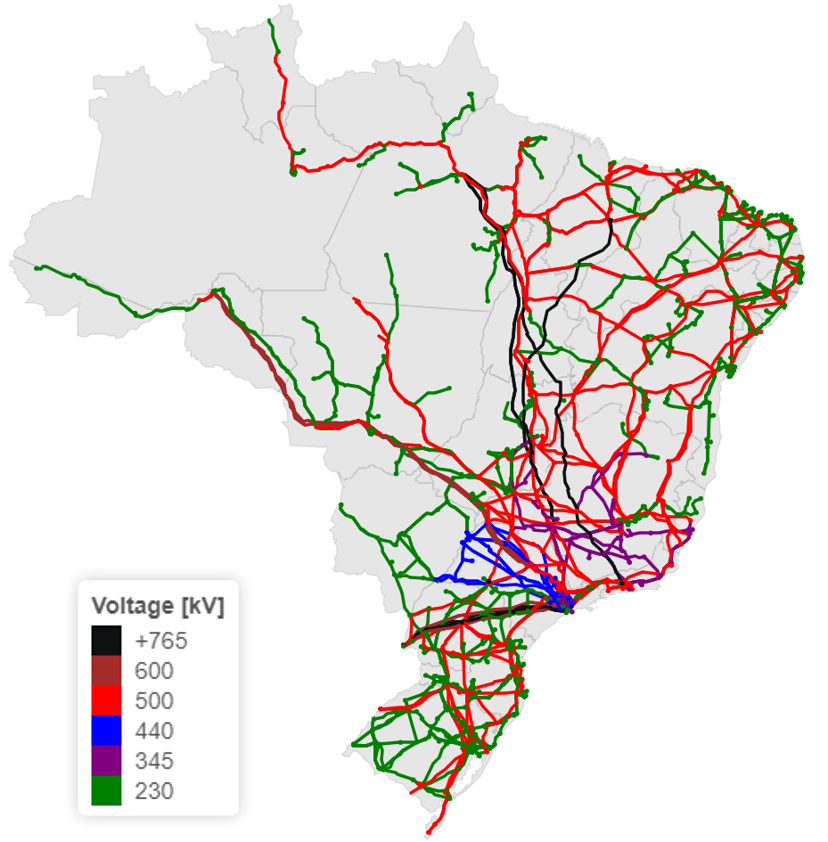}
 %where an .eps filename suffix will be assumed under Latex, 
 %and a .pdf suffix will be assumed for pdfLatex; or what has been declared
 %via \DeclareGraphicsExtensions.
\caption{Brazilian power system in the end of 2026.}
\label{fig: Brazilian network}
\end{figure}

Figure \ref{fig: Brazilian network} shows the Brazilian high voltage network of 2026. As mentioned, Brazil's grid has continental dimensions, with almost 200,000 kilometers of transmission lines, as shown in table \ref{tab: LT km 2027}. One of the main reasons for the extension of this system is the long distances between the big power plants and the load center located in the Southeast. Most of the hydroelectric potential is located in the north of the country, in the Amazon region, where are the largest hydroelectric plants, such as the Belo Monte (11,233 MW) and the Madeira River Hydro Complex, Santo Antônio (3,568 MW) and Jirau (3,750 MW). For the flow of energy from these northern hydroelectric plants, two direct current transmission systems were designed: two 600 kV DC bipoles connecting the Madeira river complex to the Araraquara substation in São Paulo; and two 800 kV DC bipoles connecting Belo Monte to Estreito and Terminal Rio substations, both in the southeastern region.

\begin{table}[!ht]
	% increase table row spacing, adjust to taste
	\renewcommand{\arraystretch}{1.3}
	\centering
	\caption{Total length of transmission lines in the end of 2026.}
	\label{tab: LT km 2027}
	\begin{tabular}{c c}
		\toprule
			Voltage [kV] & Total Length [thousands of km] \\
			\cline{1-2}
			$\pm$ 800 DC & 12.4 \\
            750 & 2.7 \\
            $\pm$ 600 DC & 12.9 \\
            500 & 76.2 \\
            440 & 7.0 \\
            345 & 11.8 \\
            230 & 73.9 \\
            Total & 196.8 \\
		\bottomrule
	\end{tabular}
\end{table}

In addition to hydroelectric generation, Brazil has a great wind source potential in the northeast and south regions and solar source in the the Northeast and Southeast. Thus, there are several 500 kV AC tie-lines among these regions to properly flow the energy from the power source to the load center. 

\subsection{Regional aggregation}
\label{sec: areas}
\begin{figure}[!ht]
	\centering
	\includegraphics[width=3in]{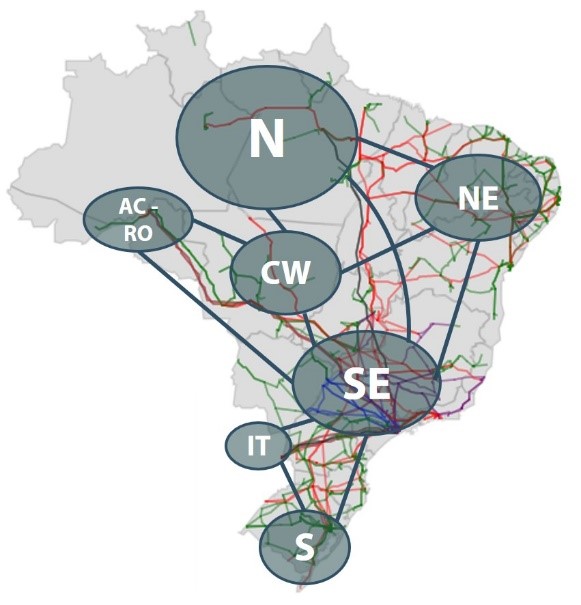}
	%where an .eps filename suffix will be assumed under Latex, 
	%and a .pdf suffix will be assumed for pdfLatex; or what has been declared
	%via \DeclareGraphicsExtensions.
	\caption{Brazilian areas}
	\label{fig: Brazilian areas}
\end{figure}

As mentioned, the transmission network was divided into seven regions (figure \ref{fig: Brazilian areas}). These regions were created by the clustering of bus marginal costs, obtained from a probabilistic simulation of the system operation. The rationale was that, if two buses have similar marginal costs, there are no significant transmission constraints between them. The clustering method described in \cite{rsampaio2019community}, also ensures that the buses in each cluster are connected.

The interchange constraints among these areas are summarized in table \ref{tab: interconnection 2027}. Long distances transmission lines often cause problems related to transient stability and reactive power compensation, limiting the power transfer capability. For this reason and reliability issues, the Brazilian system operator limits the circuit flows in the transmission network. In particular, such constraints can be observed in 500 kV circuits connecting the Northeast to the North and South-eastern regions. In fact, the ratio between the allowed flow limit and the sum of the thermal capacity of the circuits in the N/NE connection is around 30\%. It is noteworthy that, as the system gets more meshed, the exchange limits tend to increase, thus, increasing the power transfer capacity.

\begin{table}[!ht]
	% increase table row spacing, adjust to taste
	\renewcommand{\arraystretch}{1.3}
	\centering
	\caption{Interconnection capacity in MW between regions in 2026.}
	\label{tab: interconnection 2027}
	\begin{tabular}{c c c}
		\toprule
		     & From $\rightarrow$ To & To $\rightarrow$ From\\
		    \cline{2-3}
            SE - SU & 7200 & 12900 \\
            SE - NE & 6000 & 6000 \\
            SE - CO & 2825 & 28259 \\
            NO - NE & 6500 & 6500 \\
            NO - CO & 4100 & 4000 \\
            IT - SE & 5500 & 0 \\
            NO - SE & 8000 & 2500 \\
            NE - CO & 2000 & 2000 \\
            IT - SU & 8100 & 400 \\
		\bottomrule
	\end{tabular}
\end{table}

\subsection{Generation Expansion and Interconnection Reinforcements}
\label{sec: generation expansion}

As mentioned before, the first step of the planning process is to carry out a co-optimization of generation and area interconnection reinforcements. In this step, the seven-area configuration of figure \ref{fig: Brazilian areas} was used to perform the system expansion.

\begin{figure}[!ht]
	\centering
	\includegraphics[width=3.5in]{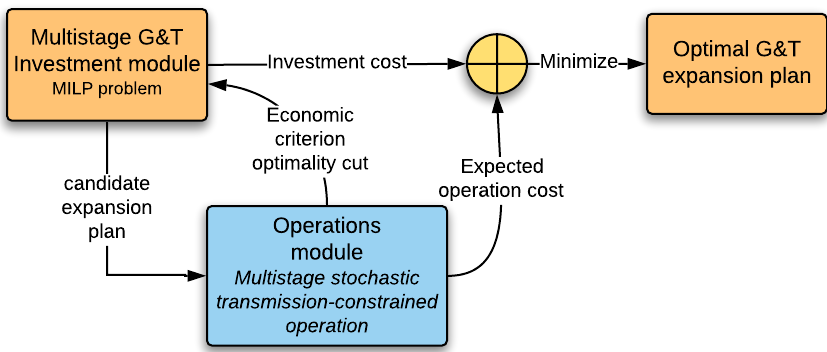}
	%where an .eps filename suffix will be assumed under Latex, 
	%and a .pdf suffix will be assumed for pdfLatex; or what has been declared
	%via \DeclareGraphicsExtensions.
	\caption{co-optmization scheme. Source \cite{fern2019stochastic}}
	\label{fig: cooptimization}
\end{figure}

The co-optimization of generation and interconnections was obtained by a Benders decomposition scheme (see figure~\ref{fig: cooptimization}), where the investment module is solved by a multistage MIP, and the operation module is solved by the SDDP multistage stochastic optimization algorithm \cite{pereira1991multi}. The probabilistic simulations were carried out in monthly steps, with three load levels in each month, for 84 hydrological and renewable generation scenarios.

\begin{figure}[!ht]
	\centering
	\includegraphics[width=3.2in]{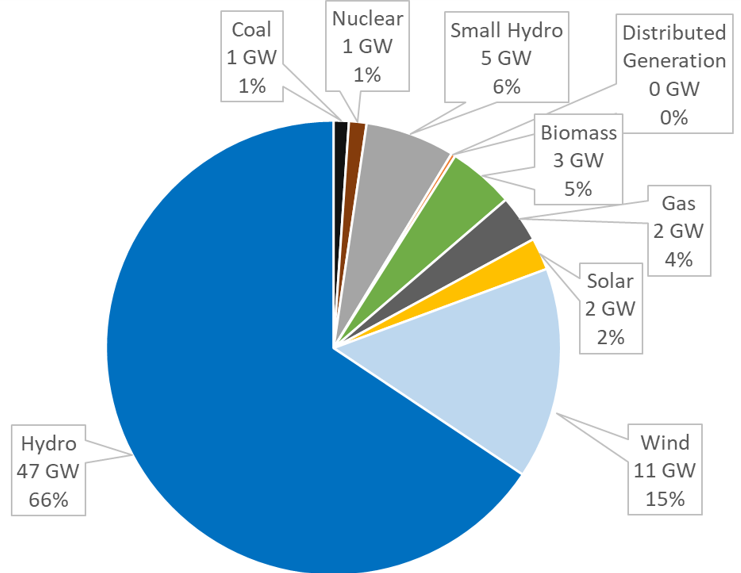}
	%where an .eps filename suffix will be assumed under Latex, 
	%and a .pdf suffix will be assumed for pdfLatex; or what has been declared
	%via \DeclareGraphicsExtensions.
	\caption{Average generation in 2027}
	\label{fig: Average generation in 2027}
\end{figure}

\begin{figure}[!ht]
	\centering
	\includegraphics[width=3.2in]{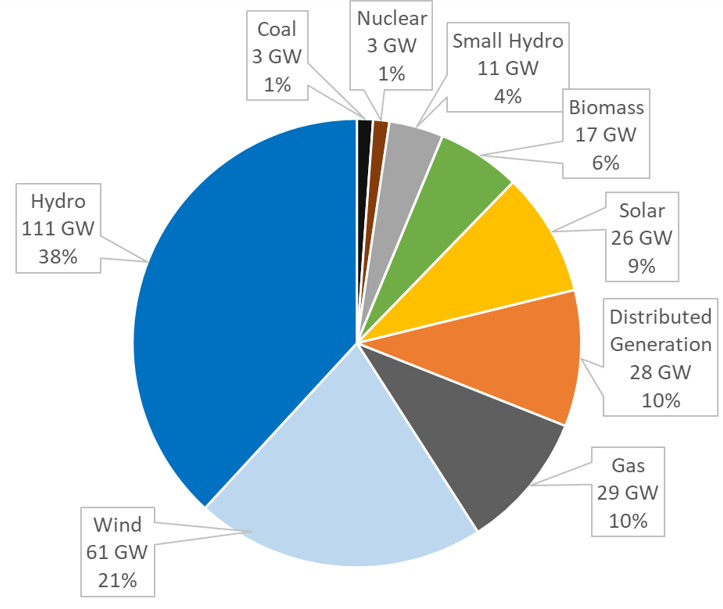}
	%where an .eps filename suffix will be assumed under Latex, 
	%and a .pdf suffix will be assumed for pdfLatex; or what has been declared
	%via \DeclareGraphicsExtensions.
	\caption{Average generation in 2036}
	\label{fig: Average generation in 2036}
\end{figure}

The initial and the final portfolio in terms of average generation selected from energy study are presented in figures \ref{fig: Average generation in 2027} and \ref{fig: Average generation in 2036}. The capacity increment by source in GW is detailed in table \ref{tab: generation}. 

\begin{table}[!ht]
	% increase table row spacing, adjust to taste
	\renewcommand{\arraystretch}{1.3}
	\centering
	\caption{Generation expansion in GW by source}
	\label{tab: generation}
    \begin{tabular}{ccccc}
    \toprule
          & \multicolumn{1}{c}{Solar} & \multicolumn{1}{c}{DG} & \multicolumn{1}{c}{Wind} & \multicolumn{1}{c}{Gas} \\
          \cline{2-5}
    2027  & \multicolumn{1}{c}{-} & 1.9  & 7.6  & \multicolumn{1}{c}{-} \\
    2028  & \multicolumn{1}{c}{-} & 2.6  & 3.1  & \multicolumn{1}{c}{-} \\
    2029  & \multicolumn{1}{c}{-} & 2.6  & 3.1  & \multicolumn{1}{c}{-} \\
    2030  & \multicolumn{1}{c}{-} & 4.0  & 1.8   & \multicolumn{1}{c}{-} \\
    2031  & 0.1  & 4.2  & 2.7  & \multicolumn{1}{c}{-} \\
    2032  & 4.3  & 3.9  & 2.5  & \multicolumn{1}{c}{-} \\
    2033  & 1.2  & 3.2  & 3.3  & \multicolumn{1}{c}{-} \\
    2034  & 1.2  & 2.4  & 5.6  & \multicolumn{1}{c}{-} \\
    2035  & 6.8  & 1.7  & 2.8  & 0.2 \\
    2036  & 2.5  & 1.1  & 3.4  & 2.2 \\
    Total & 16.1 & 27.6 & 35.9 & 2.4 \\
    \bottomrule
    \end{tabular}%
\end{table}

The leading technology in terms of capacity addition is wind power, with approximately 35 GW, mostly installed in the Northeast region, with a smaller amount in the South region. Solar photovoltaic expansion amounts almost 30 GW. In addition to the centralized solar, a total of 30 GW insertion of solar distributed generation (DG) was considered over the years, until 2036. Therefore, the total solar PV capacity addition amounts to nearly 50 GW. The generation expansion results on almost 90\% of renewable generation in 2036.

The co-optimization results of the interconnection expansion are depicted in table \ref{tab: Interconnection expansion}. Since most of the generation expansion is located in the Northeast, and the load center is in the South-east, the interchange between the Northeast and the South-east was expanded by 9GW by the end of the horizon. 

In this study, the interchange capacity increase was deployed by HVAC tie-lines and a HVDC transmission system of 4 GW, taking the massive renewable generation expansion from the Northeast straight to the load center. The HVDC system from Graça Aranha (NE) to Silvânia substation (CW) is a project that has been studied by the Brazilian Energy Research Office (EPE), and it was included do the expansion plan in the begging of the study horizon, 2027. For the HVAC tie-lines expansion, duplication of existing 500 kV transmission lines were considered in early 2032, when interconnection expansion exceeds 4 GW.

\begin{table}[!ht]
	% increase table row spacing, adjust to taste
	\renewcommand{\arraystretch}{1.3}
	\centering
	\caption{Interconnection expansion in MW from co-optimization solution.}
	\label{tab: Interconnection expansion}
	\begin{tabular}{c c}
		\toprule
		    & SE $\leftrightarrow$ NE\\
		    \cline{2-2}
			2027 & 692  \\
            2028 & 1385 \\
            2029 & 2077 \\
            2030 & 2770 \\
            2031 & 3462 \\
            2032 & 4154 \\
            2033 & 4846 \\
            2034 & 5538 \\
            2035 & 6231 \\
            2036 & 6923 \\
		\bottomrule
	\end{tabular}
\end{table}

\subsection{Transmission Expansion}
\label{sec: tranmission expansion}

Given the reinforcements of generation and area interchange capacities, the next step is to determine the optimal transmission expansion within each area. 

As mentioned,the transmission network for the year of 2026 has more than 10 thousand buses and 14 thousand circuits. In this study only the high voltage network (above 230 kV) is monitored, which is composed of 2788 circuits. 

\begin{figure}[!ht]
	\centering
	\includegraphics[width=2.7in]{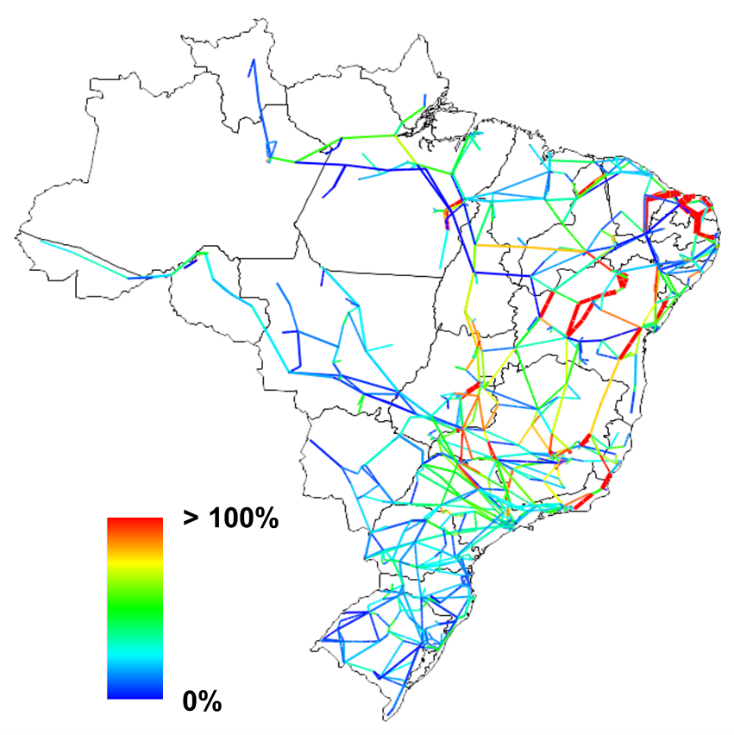}
	%where an .eps filename suffix will be assumed under Latex, 
	%and a .pdf suffix will be assumed for pdfLatex; or what has been declared
	%via \DeclareGraphicsExtensions.
	\caption{Maximum circuit loading considering all scenarios.}
	\label{fig: Max load}
\end{figure}

When solving an TEP problems one important step is an adequate selection of candidates which directly affect the MIP size. This is accomplished by analyzing the maximum loading of the circuits considering all scenarios from 2027 until 2036, as seen in figure \ref{fig: Max load}. Table~\ref{tab: num. candidates} contains the number of candidates per region. The investment cost of the transmission lines and transformers candidates were calculated from the Brazilian Electricity Regulatory Agency (ANEEL) database.

\begin{table}[!ht]
	% increase table row spacing, adjust to taste
	\renewcommand{\arraystretch}{1.3}
	\centering
	\caption{Number of circuit candidates per region.}
	\label{tab: num. candidates}
	\begin{tabular}{c c c c c c c}
		\toprule
		    S & SE & CW & NE & N & AC-RO & IT \\
			116 & 257 & 90 & 923 & 176 & 3  & 2\\
		\bottomrule
	\end{tabular}
\end{table}

The transmission expansion problem was performed using the heuristic method described in Section~\ref{sec: heuristic} from 2027 to 2036 using the predefined candidates circuits list. The candidates are the same for the entire planning horizon. However, it should be noticed that determining an adequate number of candidates per year while running a yearly TEP problem could be more suitable.

All the simulations were performed using a computer with Intel Core i7-8550U processor and 16GB RAM. The average CPU time for solving the TEP problem for a year is shown in table \ref{tab: CPU time}. 

\begin{table}[!ht]
	% increase table row spacing, adjust to taste
	\renewcommand{\arraystretch}{1.3}
	\centering
	\caption{Average CPU time .}
	\label{tab: CPU time}
	\begin{tabular}{c c c c c}
		\toprule
		    S & SE & CW & NE & N \\
			0.49 h & 1.82 h & 0.87 h & 13.93 h & 2.17 h \\
		\bottomrule
	\end{tabular}
\end{table}

\begin{figure}[!ht]
	\centering
	\includegraphics[width=3.2in]{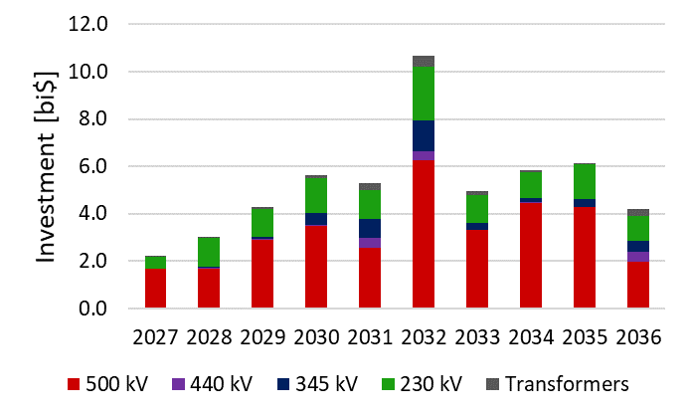}
	%where an .eps filename suffix will be assumed under Latex, 
	%and a .pdf suffix will be assumed for pdfLatex; or what has been declared
	%via \DeclareGraphicsExtensions.
	\caption{Transmission expansion by type.}
	\label{fig: expansion by type}
\end{figure}

The resulting transmission expansion plan is summarized in table \ref{tab: transmission expansion} and figure \ref{fig: expansion by type}. As expected, the Northeast received the majority of the investment due to the generation expansion investment in the area. The northern region has also been expanded, mainly to flow energy from renewable sources from the northeast to the South-east through the HVDC bipoles of Belo Monte. 

\begin{table}[!ht]
	% increase table row spacing, adjust to taste
	\renewcommand{\arraystretch}{1.3}
	\centering
	\caption{Transmission expansion decision in bi\$}
	\label{tab: transmission expansion}
    \begin{tabular}{c c c c c c c c}
    \toprule
         & S & SE & CW & NE & N & AC-RO & IT \\
        \cline{2-8}
        2027 & - & 0.05 & 0.06 & 0.96 & 1.17 & - & - \\
        2028 & 0.06 & 0.25 & 0.09 & 1.71 & 0.92 & - & - \\
        2029 & 0.06 & 0.15 & 0.14 & 1.82 & 2.11 & - & - \\
        2030 & 0.05 & 1.22 & 0.06 & 2.24 & 2.05 & - & - \\
        2031 & 0.07 & 2.03 & 0.06 & 2.18 & 0.94 & - & - \\
        2032 & 0.01 & 3.64 & 0.12 & 3.68 & 3.23 & - & - \\
        2033 & 0.03 & 0.51 & 0.09 & 2.45 & 1.88 & - & - \\
        2034 & 0.23 & 0.54 & 0.05 & 2.68 & 2.96 & - & - \\
        2035 & 0.12 & 0.65 & 0.07 & 2.16 & 3.00 & - & - \\
        2036 & 0.90 & 0.47 & 0.04 & 2.49 & 0.71 & - & - \\
    \bottomrule
    \end{tabular}%
\end{table}

Figure \ref{fig: expansion by type} shows the expansion plan depicted by transmission lines voltage level and transformers. A major investment in transmission lines of 500 kV in 2032 is noted located in the vicinity of NE/SE interconnection. This is due to the HVAC tie-lines reinforcement from the co-optimization of generation and interconnections problem solution.

\subsection{Comparison between Heuristic and Benders Decomposition Methods}
\label{sec: comparison}

Table \ref{tab: comparison} shows a comparison between the heuristic procedure and Bender decomposition solution for the first year of study (2027) for three different regions. For the Central-west area the resulting investment are exactly the same for both methods. Nevertheless, in the Southeast and the North the heuristic method resulted in a more expensive solution, since it does not guarantee the optimal solution. Yet, it is still an acceptable solution for the TEP problem, specially when analyzing large-scale network as the Brazilian power system, in which obtaining the optimal solution from Benders decomposition may not be possible due to computational performance as was the case of the Northeast region.

\begin{table}[!ht]
	% increase table row spacing, adjust to taste
	\renewcommand{\arraystretch}{1.3}
	\centering
	\caption{Comparison between Benders Decomposition and Heuristic methods results.}
	\label{tab: comparison}
    \begin{tabular}{c c c}
    \toprule
         & Heuristic Method & Benders Decomposition Method\\
         \cline{2-3}
         CO & 0.06 US\$ & 0.06 US\$ \\
         SE & 0.05 US\$ & 0.03 US\$ \\
         NO & 1.17 US\$ & 0.92 US\$ \\
         NE & 0.96 US\$ & - \\
    \bottomrule
    \end{tabular}%
\end{table}

\begin{table}[!ht]
	% increase table row spacing, adjust to taste
	\renewcommand{\arraystretch}{1.3}
	\centering
	\caption{Comparison between Benders Decomposition and Heuristic methods CPU time.}
	\label{tab: comparison time}
    \begin{tabular}{c c c}
    \toprule
          & Heuristic Method & Benders Decomposition Method\\
         \cline{2-3}
         CO & 0.87 h & 5.16 h \\
         SE & 1.82 h & 5.51 h \\
         NO & 2.17 h & 9.09 h \\
         NE & 13.93 h & time out \\
    \bottomrule
    \end{tabular}%
\end{table}

\section{Conclusions}
This paper described the new methodologies and analytical tools used to determine the least-cost expansion of Brazil’s transmission network. Planning the Brazilian transmission network expansion is a challenging analytical and computational task for several reasons, including size, investment cost and diverse operating scenarios.

First, a clustering scheme on the bus marginal costs to divide the country into a given number of areas. Then, a co-optimization of generation and interchange investment problem is solved to plan the reinforcements among those areas. This planning problem is formulated as the
co-optimization of generation and interchange investment costs, plus expected operation costs, for many hydrological and renewable scenarios. Given the planned reinforcements of generation and area interchange capacities, the next step is to determine the optimal transmission inside each area. This is accomplished by using the heuristic procedure described in section \ref{sec: heuristic}.

The effectiveness and computational performance of the proposed methodology was verified with a planning study of the Brazilian system for the period 2027-36.

\label{sec: conclusions}

% trigger a \newpage just before the given reference
% number - used to balance the columns on the last page
% adjust value as needed - may need to be readjusted if
% the document is modified later
%\IEEEtriggeratref{8}
% The 'triggered' command can be changed if desired:
%\IEEEtriggercmd{\enlargethispage{-5in}}

% references section

% can use a bibliography generated by BibTeX as a .bbl file
% BibTeX documentation can be easily obtained at:
% http://www.ctan.org/tex-archive/biblio/bibtex/contrib/doc/
% The IEEEtran BibTeX style support page is at:
% http://www.michaelshell.org/tex/ieeetran/bibtex/
%\bibliographystyle{IEEEtran}
% argument is your BibTeX string definitions and bibliography database(s)
%\bibliography{IEEEabrv,../bib/paper}
%
% <OR> manually copy in the resultant .bbl file
% set second argument of \begin to the number of references
% (used to reserve space for the reference number labels box)
%\begin{thebibliography}{1}
%\bibitem{Shell}
%M.~Shell, \emph{How to Use the IEEEtran Latex Class}, Latex Archive Contents, %\verb+http://www.ieee.org/conferences_events/+ %\verb+conferences/publishing/templates.htm+

%\bibitem{IEEEhowto:kopka}
%H.~Kopka and P.~W. Daly, \emph{A Guide to \LaTeX}, 3rd~ed.\hskip 1em plus
%  0.5em minus 0.4em\relax Harlow, England: Addison-Wesley, 1999.
  
%\end{thebibliography}

\bibliographystyle{IEEEtran}
\bibliography{myref}

% that's all folks
\end{document}